\documentclass{mcom-l}

\copyrightinfo{2010}{by the author}

\usepackage{url}

\newtheorem{theorem}{Theorem}[section]
\newtheorem{lemma}[theorem]{Lemma}
\newtheorem{proposition}[theorem]{Proposition}
\newtheorem{corollary}[theorem]{Corollary}

\theoremstyle{definition}

\newtheorem{algorithm}[theorem]{Algorithm}

\theoremstyle{remark}

\numberwithin{equation}{section}

\newcommand{\ENUMERATE}[1]
{
  \renewcommand{\labelenumi}{(\roman{enumi})}
  \begin{enumerate}#1\end{enumerate}
  \renewcommand{\labelenumi}{\arabic{enumi}.}
}
\newcommand{\IGNORE}[1]{}
\newcommand{\DEF}[1]{\stackrel{\sf def}{#1}}

\def\Q{{\mathbb Q}}

\def\F{{\mathbb F}}

\def\G{G_\alpha}

\DeclareMathOperator{\POLY}{poly}

\newcommand{\ELEMENT}[1]{\left[ #1 \right]}
\newcommand{\POWER}[2]{\ELEMENT{#1}^{#2}}
\newcommand{\LEGENDRE}[2]{\genfrac{(}{)}{}{}{#1}{#2}}
\newcommand{\nCr}[2]{\genfrac{(}{)}{0pt}{}{#1}{#2}}

\newcommand{\SET}[1]{\left\{ #1 \right\}}
\newcommand{\SETL}[2]{\left\{ #1 \;:\; #2 \right\}}
\newcommand{\SETR}[2]{\left\{ #1 \;:\; #2 \right\}}

\begin{document}

\title[On Taking Square Roots without Quadratic Nonresidues]{
On Taking Square Roots without Quadratic Nonresidues over Finite Fields
}


%    author one information
\author{Tsz-Wo Sze}
\address{Department of Computer Science, University of Maryland, College Park}
%\curraddr{}
\email{szetszwo@cs.umd.edu}
%\thanks{}

%    \subjclass is required.
%\subjclass[2000]{Primary }
%    The 2010 edition of the Mathematics Subject Classification is
%    now available.  If you are citing a classification from the
%    new scheme, use the following input coding instead.

\date{}

%\dedicatory{}

%    Abstract is required.
\begin{abstract}
We present a novel idea to compute square roots over finite fields,
without being given any quadratic nonresidue,
and without assuming any unproven hypothesis.
The algorithm is deterministic and the proof is elementary.
In some cases,
the square root algorithm runs in $\tilde{O}(\log^2 q)$ bit operations 
over finite fields with $q$ elements.
As an application,
we construct a deterministic primality proving algorithm,
which runs in $\tilde{O}(\log^3 N)$ for some integers $N$.
\end{abstract}

\maketitle

\section{Introduction}
Let $\F_q$ be a finite field with $q$ elements and $\beta\in\F_q$ be a square.
The square root problem over $\F_q$ is to find $\alpha\in\F_q$ 
such that $\alpha^2=\beta$, 
given $\beta$ and $\F_q$ as inputs.
Suppose $q\equiv 1\pmod{8}$ in this paper.
Otherwise,
the square root problem is easy;
see \cite{Crandall2001}, \cite{Cohen1993}.

The problem of taking square roots over a finite field
and the problem of constructing a quadratic nonresidue
over the same finite field
are polynomial time equivalent.
If one can take square roots,
one can compute $(-1)^{1/2}$, $(-1)^{1/4}$, $(-1)^{1/8}$, $\cdots$,
and eventually obtain a quadratic nonresidue
because the number of steps is $O(\log q)$.
Conversely, 
given a quadratic nonresidue as an additional input,
there are deterministic polynomial time algorithms
\cite{Tonelli1891}, \cite{Shanks1972} and \cite{Adleman1977}
for computing square roots.
There is no known deterministic polynomial time square root algorithm
over finite fields in general,
therefore,
there is no known deterministic polynomial time algorithm 
for constructing a quadratic nonresidue.
We discuss some probabilistic approaches below for these two problems.

There is a simple, efficient probabilistic algorithm
for finding a quadratic nonresidue
because,
in $\F_q$,
the number of quadratic nonresidues
is equal to the number of quadratic residues,
and
it is easy to determine whether an element is a quadratic nonresidue.
One could randomly pick an element $a\in\F_q$,
and then test whether $a$ is a quadratic nonresidue
by computing $a^{(q-1)/2}$.
The element $a$ is a quadratic nonresidue if and only if $a^{(q-1)/2}=-1$.
Repeat this process until a quadratic nonresidue is found.

There are several efficient probabilistic algorithms
for taking square roots in finite fields.
When quadratic nonresidues are not given,
Tonelli-Shanks \cite{Tonelli1891, Shanks1972}, 
Adleman-Manders-Miller \cite{Adleman1977}
and Cipolla-Lehmer \cite{Cipolla1903, Lehmer1969}
are considered as probabilistic algorithms 
since they require a quadratic nonresidue as an additional input.
Berlekamp-Rabin \cite{Berlekamp1970, Rabin1980}
takes square roots by polynomial factoring over finite fields.
The idea of Peralta \cite{Peralta1986}
is similar to Berlekamp-Rabin.
For other results, 
see
\cite{Bach1990},
\cite{Bach1999},
\cite{Barreto2006},
\cite{Bernstein2001},
\cite{Buchmann1996},
\cite{Muller2000},
\cite{Muller2004},
\cite{Sutherland2009}
and \cite{Turner1994}.

We restrict our discussion to prime fields $\F_p$ for odd prime $p$
in the following paragraphs.
Although there is no known deterministic polynomial time algorithm
for taking square roots,
or equivalently,
constructing a quadratic nonresidue,
over prime fields in general,
deterministic polynomial time algorithms exist in some special cases.

Schoof \cite{Schoof1985} showed 
a deterministic algorithm for computing square roots of $\beta$ over $\F_p$
with running time $O((|\beta|^{1/2+\epsilon}\log p)^9)$ bit operations\footnote
{
$|\beta|$ denotes the absolute value of $\beta$,
where $\beta$ is considered as an integer in $(-\frac{p-1}{2},\frac{p-1}{2}]$.
}
for all $\epsilon>0$.
Thus,
his algorithm is polynomial time for any fixed $\beta$
but it is exponential time generally.

For primes $p$ with $p\not\equiv 1\pmod{240}$,
a quadratic nonresidue over $\F_p$
can be constructed in deterministic polynomial time as shown below.
Denote a primitive $r$th of unity by $\zeta_r$.
If $p\not\equiv1\pmod{16}$,
at least one of
\[
\zeta_2=-1,
\qquad\zeta_4=\pm\sqrt{-1},
\qquad\zeta_8=\pm\frac{1}{\sqrt{2}}(1\pm\sqrt{-1})
\]
is a quadratic nonresidue over $\F_p$.
Suppose $p\equiv 1\pmod 4$ for the following.
If $p\equiv 2\pmod{3}$,
then the Legendre symbol $\LEGENDRE{3}{p}=\LEGENDRE{p}{3}=\LEGENDRE{2}{3}=-1$
by the law of quadratic reciprocity
and so $3$ is a quadratic nonresidue over $\F_p$.
Similarly,
$5$ is a quadratic nonresidue over $\F_p$ for $p\equiv 2,3\pmod{5}$.
Suppose $p\equiv 4\pmod{5}$.
Let
\[
\zeta_5=\frac{a+\sqrt{a^2-4}}{2},
\qquad\text{where }a=\frac{-1+\sqrt{5}}{2}\in\F_p.
\]
Then,
$a^2-4$ is a quadratic nonresidue over $\F_p$ since $\zeta_5\not\in\F_p$.
Note that the values of $\sqrt{-1}$, $\sqrt{2}$ and $\sqrt{5}$
that appeared previously
can be computed by Schoof's square root algorithm in polynomial time.
In conclusion,
the problem of constructing a quadratic nonresidue over $\F_p$ is non-trivial
only if $p\equiv 1\pmod{16}$, $p\equiv 1\pmod{3}$ and $p\equiv 1\pmod{5}$.

We end our discussion on prime fields
by considering the Extended Riemann Hypothesis (ERH).
By assuming ERH,
Ankeny \cite{Ankeny1952} showed that the least\footnote
{
The elements in $\F_p$ are considered as non-negative integers.
}
quadratic nonresidue over $\F_p$
is less than $c\log^2 p$ for some constant $c$.
As a consequence,
the probabilistic algorithm for finding a quadratic nonresidue
mentioned previously
can be improved to a deterministic polynomial time algorithm.
It can be proved that the least quadratic nonresidue must be a prime.
One could evaluate the Legendre symbol 
$\LEGENDRE{r}{p}\equiv r^{(p-1)/2}\pmod{p}$
with primes $r=2,3,5,7,...$,
until the least quadratic nonresidue is found.

In this paper,
the main results and the main ideas are presented
in Section \ref{sect-results} and Section \ref{sect-idea},
respectively.
In Section \ref{sect-group},
we construct a group and describe the arithmetic of the group.
In Section \ref{sect-sqrt},
we show a deterministic square root algorithm over finite fields.
As an application,
a deterministic primality proving algorithm is constructed
in Section \ref{sect-Primality proving}.
In the appendix (by L. Washington),
we show how to construct roots of unity needed for 
Theorem \ref{thm-sqrt 2 3^k+1}.
%%%%%%%%%%%%%%%%%%%%%%%%%%%%%%%%%%%%%%%%%%%%%%%%%%%%%%%%%%%%%%%%%%%%%%%%%%%%%%%
\section{Main Results}\label{sect-results}
We present a novel idea to compute square roots over finite fields,
without being given any quadratic nonresidue,
and without assuming any unproven hypothesis.
The square root algorithm,
Algorithm \ref{alg-sqrt},
is deterministic and the proof is elementary.
In some cases,
the algorithm runs in $\tilde{O}(\log^2 q)$ bit operations 
over finite fields $\F_q$.
As an application,
we construct a deterministic primality proving algorithm,
which runs in $\tilde{O}(\log^3 N)$ for some integers $N$.
We prove the following theorems.

\begin{theorem}\label{thm-sqrt q=2^e3^ft+1}
Let $\F_q$ be a finite field with characteristic $p$
such that 
\[
q=2^e3^ft+1
\qquad\text{and}\qquad
p\equiv 1\pmod{12}.
\]
Suppose $t=O(\POLY(\log q))$.
There is a deterministic polynomial time square root algorithm over $\F_q$.
\end{theorem}

\begin{theorem}\label{thm-sqrt 2 3^k+1} 
Let $\F_q$ be a finite field with characteristic $p$
such that 
\[
q=2^ep_1^{e_1}\cdots p_n^{e_n}t+1,
\qquad
p\equiv 13,25\pmod{36}
\qquad\text{and}\qquad
p\equiv 1\pmod{p_j},
\]
where $p_j=2\cdot3^{k_j}+1$ are $n$ distinct primes for $k_j\geq 0$.
Suppose $t+\sum p_j=O(\POLY(\log q))$.
There is a deterministic polynomial time square root algorithm over $\F_q$.
\end{theorem}

\begin{theorem}\label{thm-sqrt q=r^et+1} 
Let $\F_q$ be a finite field with characteristic $p$
such that 
\[
q=r^et+1
\qquad\text{for some prime }r.
\]
Suppose $r+t=O(\POLY(\log q))$.
There is a deterministic polynomial time square root algorithm over $\F_q$.
\end{theorem}
%%%%%%%%%%%%%%%%%%%%%%%%%%%%%%%%%%%%%%%%%%%%%%%%%%%%%%%%%%%%%%%%%%%%%%%%%%%%%%%
\section{Main Ideas}\label{sect-idea}
Suppose $\beta\in\F_q^\times$ is a square,
where $\F_q$ is a finite field with $q$ elements.
Then,
\[
\alpha^2=\beta
\qquad\text{for some }\alpha\in\F_q^\times.
\]
We present an idea to compute $\alpha$,
given $\beta$ and $\F_q$.
The problem of taking a square root of $\beta$
with arbitrary size is reduced
to the problem of constructing a primitive $r$th root of unity
$\zeta_r\in\F_q$ for some $r|q-1$.
The main ingredient of the reduction is a group isomorphism.
More details are discussed below.

Let $\G$ be a group with the following properties:
\ENUMERATE
{
\item
the group operation of $\G$ can be computed efficiently with $\beta$ 
but without the knowledge of $\alpha$,
\item
$\G$ is isomorphic to the multiplicative group $\F_q^\times$,
and
\item
the isomorphism $\psi_\alpha:\G\longrightarrow\F_q^\times$ 
depends on $\alpha$ as a parameter.
}
Since the isomorphism $\psi_\alpha$ depends on $\alpha$ 
while the value of $\alpha$ is unknown,
$\psi_\alpha$ and its inverse are not at first efficiently computable.
We try to match certain elements in $\G$
with the corresponding elements in $\F_q^\times$.
In the cases we considered,
a matched pair reveals the isomorphism $\psi_\alpha$.
Consequently,
$\alpha$ can be computed.

We first find an order $r$ element in $\G$,
where $r$ is an odd\footnote
{
The special case $r=2$ can be handled differently.
See Algorithm \ref{alg-sqrt}.
}
prime factor of $q-1$.
Write $q=r^et+1$ such that $(t,r)=1$.
Consider an element $\ELEMENT{g}\in\G$.
Suppose the order of $\ELEMENT{g}$ is $d$ such that $r|d$.
Then,
$\ELEMENT{a}=\POWER{g}{d/r}$ is an order $r$ element.
Note that there are $(r^e-1)t$ possible $\ELEMENT{g}\in\G$
leading to an order $r$  element $\ELEMENT{a}$
but only $t$ elements are not.

The element $\ELEMENT{a}$ must be matched up with
an order $r$ element in $\F_q^\times$
through the isomorphism $\psi_\alpha$.
Since $\F_q^\times$ is cyclic,
we have
\[
\psi_\alpha(\ELEMENT{a})=\zeta_r^k
\qquad
\text{for some }0<k<r,
\]
where $\zeta_r\in\F_q$ is a primitive $r$th of unity.
Once the index $k$ is obtained,
the parameter $\alpha$ of $\psi_\alpha$ can be computed.

The remaining problem is to find a primitive $r$th root of unity, $\zeta_r$.
In some cases,
$\zeta_r$ can be constructed
by taking square roots of some fixed size elements over $\F_q$.
These square roots can be computed by Schoof's square root algorithm.
In some other cases,
$\zeta_r$ can be constructed directly.

%%%%%%%%%%%%%%%%%%%%%%%%%%%%%%%%%%%%%%%%%%%%%%%%%%%%%%%%%%%%%%%%%%%%%%%%%%%%%%%
\section{A Group Isomorphic to $\F_q^\times$}\label{sect-group}
Let $\F_q$ be a finite field with $q$ odd.
Define the set
\begin{eqnarray}\label{def-G'}
\G' \DEF= \SETR{\ELEMENT{a}}{a\in\F_q,a\neq\pm\alpha}
\qquad\text{for some }\alpha\in\F_q^\times.
\end{eqnarray}
For distinguishing the elements in $\G'$ and the elements in $\F_q$,
we denote the former by $\ELEMENT{\,\cdot\,}$.
The number of elements in $\G'$ is $q-2$.
By adding the element $\ELEMENT{\infty}$ to $\G'$,
we obtain
\begin{eqnarray}\label{def-G_alpha}
\G \DEF= \G' \cup\SET{\ELEMENT{\infty}}.
\end{eqnarray}
Define an operation $*$ on $\G$ as follows:
$\forall \ELEMENT{a}\in\G$ 
and $\forall\ELEMENT{a_1},\ELEMENT{a_2}\in \G'$ 
with $a_1+a_2\neq 0$,
\begin{eqnarray}
\label{eqn-a*infty}
{\ELEMENT{a}} * \ELEMENT{\infty} &=& \ELEMENT{\infty}*\ELEMENT{a}=\ELEMENT{a},\\
\label{eqn-a*-a}
{\ELEMENT{a_1}}*\ELEMENT{-a_1} &=&\ELEMENT{\infty}, \\
\label{eqn-a_1*a_2}
{\ELEMENT{a_1}}*\ELEMENT{a_2} &=& \ELEMENT{\frac{a_1a_2+\alpha^2}{a_1+a_2}}.
\end{eqnarray}
Interestingly, 
$(\G,*)$ is a well-defined group,
which is isomorphic to the multiplicative group $\F_q^\times$.
The group $\G$ provides 
a new computational point of view of $\F_q^\times$.
We will use $\G$ to construct a deterministic square root algorithm later.

\begin{theorem}
$(\G,*)$ is an Abelian group with identity $\ELEMENT{\infty}$.
The group $\G$ is isomorphic to the multiplicative group $\F_q^\times$.
\end{theorem}
\begin{proof}
Define a bijective mapping
\begin{equation}\label{eqn-psi}
\psi: \G\longrightarrow\F_q^\times,
\qquad \ELEMENT{\infty}\longmapsto1,
\qquad \ELEMENT{a}\longmapsto\frac{a+\alpha}{a-\alpha}
\end{equation}
with inverse
\begin{equation}\label{eqn-psi^-1}
\psi^{-1}:\F_q^\times\longrightarrow \G,
\qquad 1\longmapsto\ELEMENT{\infty},
\qquad b\longmapsto\ELEMENT{\frac{\alpha(b+1)}{b-1}}.
\end{equation}
A straightforward calculation shows that $\psi$ is a homomorphism.
The Theorem follows.
\end{proof}

Note that $\G$ is cyclic because $\F_q^\times$ is.
Since $q$ is odd,
there is a unique order 2 element in $\G$.
For any $\alpha\in\F_q^\times$,
we have 
\[\psi(\ELEMENT{0})=\frac{0+\alpha}{0-\alpha}=-1.\]
Thus,
$\ELEMENT{0}$ is the order 2 element in $\G$,
independent of the choice of $\alpha$.

For more discussions on $\G$,
see \cite{stw_phd}.
%%%%%%%%%%%%%%%%%%%%%%%%%%%%%%%%%%%%%%%%%%%%%%%%%%%%%%%%%%%%%%%%%%%%%%%%%%%%%%%
\subsection{Singular Curves with a Double Point}
We can reinterpret the group law in terms of ``singular elliptic curves.''
Consider the curve
\[E: y^2=x^2(x+\alpha^2).\]
Let $E(\F_q)$ be the points on the curve with coordinates in $\F_q$.
The only singular point on $E(\F_q)$ is $(0,0)$,
which is a double point.
Let $E_{ns}(\F_q)$ be the non-singular points on $E(\F_q)$.
Then,
the mapping
\[
\tau:E_{ns}(\F_q)\rightarrow\F_q^\times,
\qquad\infty\longmapsto1,
\qquad(x,y)\longmapsto\frac{(y/x)+\alpha}{(y/x)-\alpha}
\]
is an isomorphism from $E_{ns}(\F_q)$ to $\F_q^\times$.
The inverse is
\[
\tau^{-1}:\F_q^\times\rightarrow E_{ns}(\F_q),
\qquad 1\longmapsto \infty,
\qquad\lambda\longmapsto
\left(\frac{4\alpha^2\lambda}{(\lambda-1)^2},
\frac{4\alpha^3(\lambda+1)}{(\lambda-1)^3}\right).
\]
For proofs and details, 
see \cite{lcw2008} p61 - p63.
Together with the isomorphism $\psi$ given in equation (\ref{eqn-psi}),
we have
\[
\G
\quad\simeq\quad
\F_q^\times
\quad\simeq\quad
E_{ns}(\F_q).
\]
The isomorphism from $E_{ns}(\F_q)$ to $\G$ is surprisingly simple:
\[
\psi^{-1}\circ\tau:E_{ns}(\F_q)\longrightarrow\G,
\qquad\infty\longmapsto\ELEMENT{\infty},
\qquad(x,y)\longmapsto \ELEMENT{y/x}.
\]

Although it is possible to formulate our discussion
in terms of the language of elliptic curves,
we will keep using $\G$ in this paper.
%%%%%%%%%%%%%%%%%%%%%%%%%%%%%%%%%%%%%%%%%%%%%%%%%%%%%%%%%%%%%%%%%%%%%%%%%%%%%%%
\section{Taking Square Roots}
\label{sect-sqrt}
Suppose $\beta\in\F_q^\times$ is a square.
We have
\[
\alpha^2=\beta
\qquad\text{for some }\alpha\in\F_q^\times.
\]
Consider the group $\G$ defined in equation (\ref{def-G_alpha}).
Let $\zeta_d\in\F_q$ be a primitive $d$th root of unity for $d|q-1$.
We have the following proposition.

\begin{proposition}\label{prop-alpha=a(zeta_d^k-1)/(zeta_d^k+1)}
Let $\ELEMENT{a}\in\G$ such that $\ELEMENT{a}^2\neq\ELEMENT{\infty}$.
Suppose $\ELEMENT{a}^d=\ELEMENT{\infty}$ for some $d>0$.
Then,
\[
\alpha=\pm\frac{a(\zeta_d^k-1)}{\zeta_d^k+1}
\qquad\text{for some }0<k<\frac{d}{2}.
\]
\end{proposition}
\begin{proof}
Since $\psi$ defined in equation (\ref{eqn-psi}) is an isomorphism,
we have
\[
\psi(\ELEMENT{a})^d
=\psi(\ELEMENT{a}^d)
=\psi(\ELEMENT{\infty})
=1
\]
over $\F_q^\times$.
Then
\[
\psi(\ELEMENT{a})=\zeta_d^j
\qquad\text{for some }0\leq j<d.
\]
Since $\ELEMENT{a}^2\neq\ELEMENT{\infty}$ by assumption,
we have $j\neq 0$ and $j\neq\frac{d}{2}$.
By applying $\psi^{-1}$ on both sides,
we obtain
\[
\ELEMENT{a}=\psi^{-1}(\zeta_d^j)
=\ELEMENT{\alpha(\zeta_d^j+1)/(\zeta_d^j-1)}.
\]
Therefore,
\[
\alpha=a(\zeta_d^j-1)/(\zeta_d^j+1).
\]
If $j<\frac{d}{2}$,
the proposition follows by setting $k=j$.
If $j>\frac{d}{2}$,
let $k=d-j<\frac{d}{2}$.
Then,
\[\frac{a(\zeta_d^k-1)}{\zeta_d^k+1}
=\frac{a(\zeta_d^{-j}-1)}{\zeta_d^{-j}+1}
=\frac{a(1-\zeta_d^j)}{1+\zeta_d^j}=-\alpha.\]
The Proposition follows.
\end{proof}

Proposition \ref{prop-alpha=a(zeta_d^k-1)/(zeta_d^k+1)}
suggests a method to compute $\alpha$.
The ingredients are
(1) an element $\ELEMENT{a}\in\G$ such that $\ELEMENT{a}^d=\ELEMENT{\infty}$,
(2) a primitive $d$th root of unity $\zeta_d\in\F_q$,
and
(3) the index $k$.
It also requires that
the power $\ELEMENT{a}^k$ has to be efficiently computable.
Recall that $\G'$,
which is defined in equation (\ref{def-G'}),
is the set of all elements in $\G$ except the identity.

\begin{lemma}\label{lem-computing [g1]*[g2]}
Given a square $\beta\in\F_q^\times$,
the group operation $*$ over $\G$ can be performed
in $\tilde{O}(\log q)$ bit operations
without the knowledge of $\alpha$.
\end{lemma}
\begin{proof}
Clearly,
the group operation involving the identity element is trivial.
By equations (\ref{eqn-a*-a}) and (\ref{eqn-a_1*a_2}),
for any $\ELEMENT{g_1},\ELEMENT{g_2}\in\G'$,
\begin{equation}\label{eqn-g_1*g_2}
\ELEMENT{g_1}*\ELEMENT{g_2}=
\begin{cases}
\ELEMENT{\infty} &\text{, if }g_1+g_2=0, \\
\ELEMENT{\frac{g_1g_2+\beta}{g_1+g_2}} &\text{, otherwise.}
\end{cases}
\end{equation}
Note that equation (\ref{eqn-g_1*g_2}) does not involve $\alpha$.
Therefore,
the group operation $*$ over $\G$ can be computed
by a few field operations over $\F_q$ in the worst case.
The Lemma follows from the fact that
field operations over $\F_q$ can be performed in $\tilde{O}(\log q)$ bit operations;
see \cite{Furer2007}, \cite{Schoenhage1971}, \cite{Knuth1969},
\cite{Gathen2003}.
\end{proof}

\begin{lemma}\label{lem-computing [g]^k}
Given a square $\beta\in\F_q^\times$,
the power $\POWER{g}{k}$ for any $\ELEMENT{g}\in\G'$
can be computed in $\tilde{O}(\log k\log q)$ bit operations
without the knowledge of $\alpha$.
\end{lemma}
\begin{proof}
The power $\ELEMENT{g}^k$ can be evaluated in $O(\log k)$ group operations
using the successive squaring method.
The Lemma follows from Lemma \ref{lem-computing [g1]*[g2]}.
\end{proof}

\IGNORE
{
Another way to prove the lemma
is due to Proposition \ref{prop-[a]^d=infty iff Psi_d(a)=0}
and equation (\ref{eqn-[a]^k}).
If $\Psi_k(g)=0$,
then $\ELEMENT{g}^k=\ELEMENT{\infty}$.
Otherwise,
$\ELEMENT{g}^k=\ELEMENT{\frac{\gamma_k(g)}{\Psi_k(g)}}$.
The values of $\gamma_k(g)$ and $\Psi_k(g)$
can be evaluated by the recursion equations 
in Lemma \ref{lem-gamma,Psi k+1 recursion}
and Lemma \ref{lem-gamma,Psi 2k recursion}.
Note that the coefficients in $\gamma$'s, $\Psi$'s and the recursion equations
only involve integers and $\beta$, 
but not $\alpha$.
The running time is also $\tilde{O}(\log k\log q)$.
}
%%%%%%%%%%%%%%%%%%%%%%%%%%%%%%%%%%%%%%%%%%%%%%%%%%%%%%%%%%%%%%%%%%%%%%%%%%%%%%%
\subsection{The Algorithm}\label{sect-the algorithm}
In this section,
we present a deterministic square root algorithm over $\F_q$.
Write
\begin{equation}\label{eqn-q}
q=2^{e}p_1^{e_1}\cdots p_n^{e_n}t+1,
\end{equation}
where $p_1,\cdots,p_n$ are $n$ distinct odd primes
and $t,e,e_1,\cdots,e_n$ are positive integers
such that $(2p_1\cdots p_n,t)=1$.
Suppose $e>1$.
Otherwise,
the square root problem is easy.
We have following algorithm.

\begin{algorithm}[Taking Square Roots]\label{alg-sqrt}
The inputs are $\beta$ and $\F_q$,
where $\beta\in\F_q^\times$ is a square.
This algorithm returns $\pm\sqrt{\beta}$.
\end{algorithm}
\begin{enumerate}
\item[I.]
Consider $2t-1$ distinct elements $g_1,g_2,...,g_{2t-1}\in\F_q^\times$.
\begin{enumerate}
\item[I.1]
If there exists $t'$ such that $g_{t'}^2=\beta$,
then return $\pm g_{t'}$.
\item[I.2]
Otherwise,
set $g=g_{t''}$ for some $t''$
such that $\POWER{g_{t''}}{2t}\neq\ELEMENT{\infty}$.
\end{enumerate}
\item[II.]
If $\POWER{g}{(q-1)/2^{e-1}}\neq\ELEMENT{\infty}$,
do the following:
\begin{enumerate}
\item[II.1]
Find the largest $k$ such that 
$\POWER{g}{(q-1)/2^k}=\ELEMENT{\infty}$.
\item[II.2]
Compute $\ELEMENT{a}=\POWER{g}{(q-1)/2^{k+2}}$,
an order $4$ element in $\G$.
\item[II.3]
Return $\pm a\sqrt{-1}$.
\end{enumerate}
\item[III.]
Find $m$ such that $\POWER{g}{(q-1)/p_{m}^{e_{m}}}\neq\ELEMENT{\infty}$.
\item[IV.]
Set $r=p_{m}$ and then do the following:
\begin{enumerate}
\item[IV.1]
Find the largest $k$
such that $\POWER{g}{(q-1)/r^k}=\ELEMENT{\infty}$.
\item[IV.2]
Compute $\ELEMENT{a}=\POWER{g}{(q-1)/r^{k+1}}$,
an order $r$ element in $\G$.
\item[IV.3]
Compute $\zeta=\zeta_r\in\F_q$, a primitive $r$th root of unity.
\item[IV.4]
Find $j$
such that $\left(a(\zeta^j-1)/(\zeta^j+1)\right)^2=\beta$
for $1\leq j\leq\frac{r-1}{2}$.
\item[IV.5]
Return $\pm a(\zeta^j-1)/(\zeta^j+1)$.
\end{enumerate}
\end{enumerate}

\begin{theorem}
Algorithm \ref{alg-sqrt} returns the square roots of $\beta$.
\end{theorem}
\begin{proof}
Clearly,
if $t'$ exists in Step I.1,
the algorithm returns the square roots of $\beta$.
Otherwise,
$\ELEMENT{g_1},\ELEMENT{g_2},\cdots,\ELEMENT{g_{2t-1}}$
are elements in $\G$.
There are $2t+1$ distinct elements
\[
\ELEMENT{\infty},\ELEMENT{0},
\ELEMENT{g_1},\ELEMENT{g_2},\cdots,\ELEMENT{g_{2t-1}}\in\G.
\]
Since $\G$ is cyclic,
the $2t$-torsion subgroup
\[
H\DEF{=}\SETR{\ELEMENT{a}\in\G}{\POWER{a}{2t}=\ELEMENT{\infty}}
\]
has exactly $2t$ elements.
We have $\ELEMENT{\infty},\ELEMENT{0}\in H$.
Therefore,
there exists $t''$ such that $\ELEMENT{g_{t''}}\not\in H$.
In Step I.2,
we obtain $g=g_{t''}$ such that $\ELEMENT{g}\not\in H$.
Denote the order of $\ELEMENT{g}$ by $d$ for the rest of the proof.

In Step II,
if $\ELEMENT{g}^{(q-1)/2^{e-1}}\neq \ELEMENT{\infty}$,
there exists $0\leq k<e-1$ 
such that 
\[
\ELEMENT{g}^{(q-1)/2^k}=\ELEMENT{\infty}
\qquad\text{and}\qquad
\ELEMENT{g}^{(q-1)/2^{k+1}}\neq\ELEMENT{\infty}.
\]
In Step II.2,
the order of $\ELEMENT{a}=\ELEMENT{g}^{(q-1)/2^{k+2}}\in\G$ is 4.
The algorithm returns 
\[
\pm a\sqrt{-1}
=\pm\frac{a(\zeta_4-1)}{\zeta_4+1},
\]
which are the square roots of $\beta$
by Proposition \ref{prop-alpha=a(zeta_d^k-1)/(zeta_d^k+1)}.
If $\ELEMENT{g}^{(q-1)/2^{e-1}}=\ELEMENT{\infty}$,
we have 
\begin{equation}\label{eqn-d|(q-1)/2^(e-1)}
d|2tp_1^{e_1}\cdots p_n^{e_n}.
\end{equation}

In Step III,
such $m$ exists.
Otherwise,
suppose $\POWER{g}{(q-1)/p_{m}^{e_{m}}}=\ELEMENT{\infty}$ for all $m$.
Then 
\[
d|(q-1)/p_{m}^{e_{m}}
\qquad\text{for all }m.
\]
Hence,
\[
d|2^et.
\]
Together with (\ref{eqn-d|(q-1)/2^(e-1)}),
we have
\[
d|2t,
\]
which contradicts $\ELEMENT{g}\not\in H$.

Step IV is similar to Step II.
Since $\ELEMENT{g}^{(q-1)/r^{e_{m}}}\neq \ELEMENT{\infty}$,
there exists $0\leq k<e_{m}$ 
such that 
\[
\ELEMENT{g}^{(q-1)/r^k}=\ELEMENT{\infty}
\qquad\text{and}\qquad
\ELEMENT{g}^{(q-1)/r^{k+1}}\neq\ELEMENT{\infty}.
\]
The order of $\ELEMENT{a}=\ELEMENT{g}^{(q-1)/r^{k+1}}\in\G$ is $r$
in Step IV.2.
By Proposition \ref{prop-alpha=a(zeta_d^k-1)/(zeta_d^k+1)},
\[
\alpha=\pm\frac{a(\zeta^j-1)}{\zeta^j+1}
\qquad\text{for some }1\leq j\leq\frac{r-1}{2}.
\]

The Theorem follows.
\end{proof}

\begin{proposition}\label{prop-sqrt-running time}
Algorithm \ref{alg-sqrt} runs in
\[
\tilde{O}((t\log t + p_{\max} + n\log q)\log q + Z_{\max})
\]
bit operations,
where $p_{\max}=\max(p_1,...,p_n)$
and $Z_{\max}=\max(Z_4,Z_{p_1},...,Z_{p_n})$,
where $Z_d$ is the time required
to construct a $d$th root of unity over $\F_q$.
\end{proposition}
\begin{proof}
Writing $q$ in the form of equation (\ref{eqn-q}) by trial divisions
requires $\tilde{O}(p_{\max}\log q)$.

The running time of Step I is $\tilde{O}(t\log t\log q)$
since multiplications over $\F_q$ and powering over $\G$ 
can be performed in $\tilde{O}(\log q)$ and $\tilde{O}(\log t\log q)$,
respectively.

In Step II,
computing $\POWER{g}{(q-1)/2^{e-1}}$,
finding the required $k$ in Step II.1 and computing $\ELEMENT{a}$ in Step II.2
take $\tilde{O}(\log^2 q)$.
It also requires $O(Z_4)$ to compute $\zeta_4=\sqrt{-1}$.
The running time of Step II is $\tilde{O}(\log^2 q+Z_4)$.

Clearly,
the running time of Step III is $\tilde{O}(n\log^2 q)$.

Step IV is similar to Step II
except that there are $(r-1)/2$ possible $j$ in Step IV.4,
which takes $\tilde{O}(r\log q)$.
The running time of Step IV is $\tilde{O}((r+\log q)\log q + Z_r)$.

The Proposition follows.
\end{proof}

\begin{corollary}\label{cor-sqrt-polynomial time}
Algorithm \ref{alg-sqrt} runs in polynomial time
when 
\[
t+p_{\max}+Z_{\max}=O(\POLY(\log q)).
\]
\end{corollary}
\begin{proof}
This immediately follows from Proposition \ref{prop-sqrt-running time}.
\end{proof}

We consider some special cases in the rest of the section.
%%%%%%%%%%%%%%%%%%%%%%%%%%%%%%%%%%%%%%%%%%%%%%%%%%%%%%%%%%%%%%%%%%%%%%%%%%%%%%%
\subsection{Case $q=2^e3^ft+1$}
Consider the finite fields $\F_q$ with characteristic $p$
such that $q=2^e3^ft+1$ and $p\equiv 1\pmod{12}$.
Note that $e\geq 2$ and $f\geq 1$ because $p\equiv 1\pmod{12}$.
We prove Theorem \ref{thm-sqrt q=2^e3^ft+1} below.

\begin{proof}[Proof of Theorem \ref{thm-sqrt q=2^e3^ft+1}]
The elements $-1$ and $-3$ are squares in the prime field $\F_p$.
We can compute $\zeta_3=\frac{-1\pm\sqrt{-3}}{2}$ and $\zeta_4=\sqrt{-1}$
in $\tilde{O}(\log^9p)$ by Schoof's square root algorithm.
Then, 
the running time of Algorithm \ref{alg-sqrt} is
\[
\tilde{O}((t\log t+\log q)\log q + \log^9p)
\]
bit operations
by Proposition \ref{prop-sqrt-running time}.
Since $t=O(\POLY(\log q))$ by assumption,
the Theorem follows.
\end{proof}
%%%%%%%%%%%%%%%%%%%%%%%%%%%%%%%%%%%%%%%%%%%%%%%%%%%%%%%%%%%%%%%%%%%%%%%%%%%%%%%
\subsection{Constructing Primitive $(2\cdot3^k+1)$th Roots of Unity}
Suppose $p$ be a prime with $p\equiv1\pmod{4}$ and $p\equiv 4,7\pmod{9}$.
We show in Lemma \ref{lem-cube roots} below that
cube roots over $\F_p$ can be computed efficiently.
As a consequence,
a primitive $r$th root of unity $\zeta_r$,
for prime $r=2\cdot3^k+1$ and some $k\geq 1$,
can be computed in polynomial time by the method described in the Appendix.
We will prove Theorem \ref{thm-sqrt 2 3^k+1} after Lemma \ref{lem-cube roots}.

\begin{lemma}\label{lem-cube roots}
Let $p$ be a prime with $p\equiv1\pmod{4}$ and $p\equiv 4,7\pmod{9}$.
Cube roots over $\F_p$ can be computed in polynomial time.
\end{lemma}
\begin{proof}
We can compute $\zeta_3=\frac{-1\pm\sqrt{-3}}{2}\in\F_p$
by Schoof's square root algorithm.
Let $b\in\F_p$ be a cubic residue.
We have $b^{(p-1)/3}=1$.
If $p\equiv 4\pmod{9}$,
let $a=b^{(2p+1)/9}$.
Then,
\[
a^3=b^{(2p+1)/3}=b^{1+2(p-1)/3}=b.
\]
Therefore,
$b^{(2p+1)/9}$, $b^{(2p+1)/9}\zeta_3$ and $b^{(2p+1)/9}\zeta_3^2$
are cube roots of $b$.
Similarly,
if $p\equiv 7\pmod{9}$,
let $a=b^{(p+2)/9}$.
Then, 
\[
a^3=b^{(p+2)/3}=b^{1+(p-1)/3}=b.
\]
Therefore,
$b^{(p+2)/9}$, $b^{(p+2)/9}\zeta_3$ and $b^{(p+2)/9}\zeta_3^2$
are cube roots of $b$.
All computations can be performed in polynomial time.
The lemma follows.
\end{proof}

\begin{proof}[Proof of Theorem \ref{thm-sqrt 2 3^k+1}]
The square roots $\sqrt{-1},\sqrt{p_1},\sqrt{p_2},...,\sqrt{p_n}\in\F_p$
can be computed using Schoof's square root algorithm.
Since $p\equiv 13,25\pmod{36}$,
cube roots over $\F_p$ can be computed in polynomial time 
by Lemma \ref{lem-cube roots}.
Then,
the primitive roots of unity
$\zeta_{p_1},\zeta_{p_2},...,\zeta_{p_n}\in\F_p$ can be constructed
by the method described in the Appendix.
Since $t+\sum p_j=O(\POLY(\log q))$ by assumption,
Algorithm \ref{alg-sqrt} runs in polynomial time
by Corollary \ref{cor-sqrt-polynomial time}.
The Theorem follows.
\end{proof}
%%%%%%%%%%%%%%%%%%%%%%%%%%%%%%%%%%%%%%%%%%%%%%%%%%%%%%%%%%%%%%%%%%%%%%%%%%%%%%%
\subsection{Searching for Primitive Roots of Unity}\label{sect-searching zeta_r}
In the previous sections,
the square root problem with arbitrary size elements is first reduced to
the problem of constructing primitive roots of unity,
which is further reduced to the square root problem
with some fixed size elements.
We show in Algorithm \ref{alg-zeta} below that
a primitive root of unity can be constructed efficiently
without the need of taking square roots in some cases.
We will prove Theorem \ref{thm-sqrt q=r^et+1} at the end of the section.

\begin{algorithm}[Constructing a Primitive $r$th Root of Unity]
\label{alg-zeta}
The inputs are $r$ and $\F_q$
for some odd prime $r$
such that $q=r^et+1$ and $(r,t)=1$.
This algorithm returns a primitive $r$th root of unity in $\F_q$.
\end{algorithm}
\begin{enumerate}
\item[1.]
Consider $t+1$ distinct elements $g_1,...,g_{t+1}\in\F_q^\times$. \\
Set $g=g_j$ such that $g_j^t\neq1$.
\item[2.]
Find the largest $k$ such that $g^{(q-1)/r^k}=1$.
\item[3.]
Return $g^{(q-1)/r^{k+1}}$.
\end{enumerate}

\begin{lemma}\label{lem-zeta-correctness}
Algorithm \ref{alg-zeta} returns a primitive $r$th root of unity.
\end{lemma}
\begin{proof}
Since the $t$-torsion subgroup of $\F_q^\times$ only has $t$ elements
but there are $t+1$ distinct elements in $g_1,g_2,\cdots,g_{t+1}$,
there exists an element $g_j$ such that $g_j^t\neq 1$.
Let $d$ be the order of $g=g_j$.
Then,
$r$ divides $d$ 
and there exists $k$ such that 
\[
g^{(q-1)/r^k}=1
\qquad{and}\qquad
g^{(q-1)/r^{k+1}}\neq 1,
\]
which means that $g^{(q-1)/r^{k+1}}$ is a primitive $r$th root of unity.
\end{proof}

\begin{lemma}\label{lem-zeta-running time}
Algorithm \ref{alg-zeta} runs in $\tilde{O}((t\log t+\log q)\log q)$ bit operations.
\end{lemma}
\begin{proof}
The running time for Step 1 is $\tilde{O}(t\log t\log q)$
and the running time for Step 2 and Step 3 is $\tilde{O}(\log^2 q)$.
The Lemma follows.
\end{proof}

Similarly,
we may construct a primitive $4$th root of unity
by Algorithm \ref{alg-zeta_4} below.
The correctness proof for Algorithm \ref{alg-zeta_4}
is similar to the proof given for Lemma \ref{lem-zeta-correctness}.
The running time is also $\tilde{O}((t\log t+\log q)\log q)$.

\begin{algorithm}[Constructing a Primitive $4$th Root of Unity]
\label{alg-zeta_4}
The input is $\F_q$
such that $q=2^et+1$,
where $e>1$ and $t$ is odd.
This algorithm returns a primitive $4$th root of unity in $\F_q$.
\end{algorithm}
\begin{enumerate}
\item[1.]
Consider $2t+1$ distinct elements $g_1,...,g_{2t+1}\in\F_q^\times$. \\
Set $g=g_j$ such that $g_j^{2t}\neq1$.
\item[2.]
Find the largest $k$ such that $g^{(q-1)/2^k}=1$.
\item[3.]
Return $g^{(q-1)/2^{k+2}}$.
\end{enumerate}

\begin{proof}[Proof of Theorem \ref{thm-sqrt q=r^et+1}]
If $r=2$,
construct $\zeta_4$ by Algorithm \ref{alg-zeta_4}.
Otherwise,
construct $\zeta_r$ by Algorithm \ref{alg-zeta}.
The running time is $\tilde{O}((t\log t+\log q)\log q)$ for either case.
Then,
the running time of Algorithm \ref{alg-sqrt}
is
\[
\tilde{O}((t\log t+r+\log q)\log q)
\]
bit operations
by Proposition \ref{prop-sqrt-running time}.
Since $r+t=O(\POLY(\log q))$ by assumption,
the Theorem follows.
\end{proof}

\section{Deterministic Primality Proving}\label{sect-Primality proving}
We briefly describe a deterministic primality proving algorithm
as an application of the square root algorithm.
For more details,
see \cite{stw_proth}.

Suppose $N=2^et+1>3$ for some odd $t$ with $2^e>t$.
Try to compute $\sqrt{-1}$ by Algorithm \ref{alg-zeta_4}
and then try to compute
$(-1)^{1/4}$, $(-1)^{1/8}$, $\cdots$, $(-1)^{1/2^{e-1}}$
by Algorithm \ref{alg-sqrt}.
If $(-1)^{1/2^{e-1}}$ is obtained,
then $N$ is a prime by Proth's Theorem (Theorem \ref{thm-Proth's} below).
Otherwise,
since the square root algorithm is deterministic,
the computation process must fail in some point
and then we conclude that $N$ is composite.
Such a primality proving algorithm is deterministic 
and runs in 
\[
\tilde{O}((t\log t + \log N)\log^2 N)
\]
bit operations.

The algorithm runs in $\tilde{O}(\log^3 N)$
when $t$ is $O(\log N)$.
For this kind of numbers,
the algorithm is faster than other applicable deterministic algorithms.
The running time of the AKS algorithm \cite{Agrawal2004} 
and Lenstra-Pomerance's modified AKS algorithm \cite{Lenstra2009}
are $\tilde{O}(\log^{7.5} N)$ and $\tilde{O}(\log^6 N)$, 
respectively.
Assuming ERH,
Miller's algorithm \cite{Miller1975} is deterministic 
with running time $\tilde{O}(\log^4N)$.

\begin{theorem}\label{thm-Proth's}
{\bf (Proth's Theorem)}
Let $N=2^et+1$ for some odd $t$ with $2^e>t$.
If
\[a^{(N-1)/2}\equiv -1\pmod{N}\]
for some $a$,
then $N$ is a prime.
\end{theorem}
See \cite{Williams1998} for the details of Proth's Theorem.
\bibliographystyle{amsplain}
\bibliography{sqrt}

\providecommand{\bysame}{\leavevmode\hbox to3em{\hrulefill}\thinspace}
\providecommand{\MR}{\relax\ifhmode\unskip\space\fi MR }
\providecommand{\MRhref}[2]{%
  \href{http://www.ams.org/mathscinet-getitem?mr=#1}{#2}
}
\providecommand{\href}[2]{#2}
\begin{thebibliography}{10}

\bibitem{Adleman1977}
Leonard~M. Adleman, Kenneth~L. Manders, and Gary~L. Miller, \emph{On taking
  roots in finite fields}, Proceedings of the 18th IEEE Symposium on
  Foundations of Computer Science, IEEE, 1977, pp.~175--178.

\bibitem{Agrawal2004}
Manindra Agrawal, Neeraj Kayal, and Nitin Saxena, \emph{{PRIMES} is in {P}},
  Ann. of Math. \textbf{160} (2004), no.~2, 781--793.

\bibitem{Ankeny1952}
Nesmith~C. Ankeny, \emph{The least quadratic non residue}, Ann. of Math.
  \textbf{55} (1952), no.~1, 65--72.

\bibitem{Bach1990}
Eric Bach, \emph{A note on square roots in finite fields}, IEEE Transactions on
  Information Theory \textbf{36} (1990), no.~6, 1494--1498.

\bibitem{Bach1999}
Eric Bach and Klaus Huber, \emph{Note on taking square-roots modulo {$N$}},
  IEEE Transactions on Information Theory \textbf{45} (1999), no.~2, 807--809.

\bibitem{Barreto2006}
Paulo S. L.~M. Barreto and Jos{\'e}~Felipe Voloch, \emph{Efficient computation
  of roots in finite fields}, Des. Codes Cryptography \textbf{39} (2006),
  no.~2, 275--280.

\bibitem{Berlekamp1970}
Elwyn~R. Berlekamp, \emph{Factoring polynomials over large finite fields},
  Math. Comp. \textbf{24} (1970), no.~111, 713--735.

\bibitem{Bernstein2001}
Daniel~J. Bernstein, \emph{Faster square roots in annoying finite fields},
  2001, preprint (\url{http://cr.yp.to/papers/sqroot.pdf}).

\bibitem{Buchmann1996}
Johannes Buchmann and Victor Shoup, \emph{Constructing nonresidues in finite
  fields and the extended {Riemann} hypothesis}, Math. Comp. \textbf{65}
  (1996), no.~215, 1311--1326.

\bibitem{Cipolla1903}
Michele Cipolla, \emph{Un metodo per la risoluzione della congruenza di secondo
  grado}, Napoli Rend. \textbf{9} (1903), 154--163.

\bibitem{Cohen1993}
Henri Cohen, \emph{{A Course in Computational Algebraic Number Theory}},
  Springer-Verlag, Berlin, 1993.

\bibitem{Crandall2001}
Richard Crandall and Carl Pomerance, \emph{{Prime Numbers: A Computational
  Perspective}}, Springer-Verlag, New York, 2001.

\bibitem{Furer2007}
Martin F{\"u}rer, \emph{Faster integer multiplication}, Proceedings of the 39th
  Annual ACM Symposium on Theory of Computing, ACM, 2007, pp.~57--66.

\bibitem{Lenstra2009}
Hendrik W.~Lenstra Jr. and Carl Pomerance, \emph{Primality testing with
  {Gaussian} periods}, 2009, preprint
  (\url{http://math.dartmouth.edu/~carlp/aks102309.pdf}).

\bibitem{Knuth1969}
Donald~E. Knuth, \emph{{The Art of Computer Programming, Volume 2:
  Seminumerical Algorithms}}, Addison-Wesley, Reading, 1969.

\bibitem{Lehmer1969}
Derrick~H. Lehmer, \emph{Computer technology applied to the theory of numbers},
  Studies in number theory (Englewood Cliffs, New Jersey) (William~J. Leveque,
  ed.), MAA Studies in Mathematics, vol.~6, Prentice-Hall, 1969, pp.~117--151.

\bibitem{Miller1975}
Gary~L. Miller, \emph{Riemann's hypothesis and tests for primality},
  Proceedings of Seventh Annual Symposium on Theory of Computing, ACM, 1975,
  pp.~234--239.

\bibitem{Muller2000}
Siguna M{\"u}ller, \emph{On probable prime testing and the computation of
  square roots mod n}, Algorithmic Number Theory, 4th International Symposium,
  ANTS-IV, Lecture Notes in Computer Science, vol. 1838, Springer-Verlag, 2000,
  pp.~423--437.

\bibitem{Muller2004}
\bysame, \emph{On the computation of square roots in finite fields}, Des. Codes
  Cryptography \textbf{31} (2004), no.~3, 301--312.

\bibitem{Rabin1980}
Michael~O. Rabin, \emph{Probabilistic algorithms in finite fields}, SIAM J.
  Comput. \textbf{9} (1980), no.~2, 273--280.

\bibitem{Peralta1986}
{Ren\'{e} C. Peralta}, \emph{A simple and fast probabilistic algorithm for
  computing square roots modulo a prime number}, IEEE Transactions on
  Information Theory \textbf{32} (1986), no.~6, 846--847.

\bibitem{Schoenhage1971}
Arnold Sch{\"o}nhage and Volker Strassen, \emph{{Schnelle Multiplikation gro\ss
  er Zahlen}}, Computing \textbf{7} (1971), 281--292.

\bibitem{Schoof1985}
Ren\'{e} Schoof, \emph{Elliptic curves over finite fields and the computation
  of square roots $\operatorname{mod} p$}, Math. Comp. \textbf{44} (1985),
  no.~170, 483--494.

\bibitem{Shanks1972}
Daniel Shanks, \emph{Five number-theoretic algorithms}, Proc. 2nd Manitoba
  Conf. Numer. Math., 1972, pp.~51--70.

\bibitem{Sutherland2009}
Andrew~V. Sutherland, \emph{Structure computation and discrete logarithms in
  finite abelian $p$-groups}, 2009, preprint
  (\url{http://arxiv.org/abs/0809.3413}).

\bibitem{stw_phd}
Tsz-Wo Sze, \emph{On solving univariate polynomial equations over finite fields
  and some related problems}, Ph.D. thesis, University of Maryland, 2007.

\bibitem{stw_proth}
\bysame, \emph{Deterministic primality proving on {Proth} numbers}, 2010,
  preprint (\url{http://arxiv.org/abs/0812.2596}).

\bibitem{Tonelli1891}
Alberto Tonelli, \emph{{Bemerkung \"uber die Aufl\"osung quadratischer
  Congruenzen}}, Nachrichten der Akademie der Wissenschaften in G\"ottingen
  (1891), 344--346.

\bibitem{Turner1994}
Stephen~M. Turner, \emph{Square roots $\operatorname{mod} p$}, The American
  Mathematical Monthly \textbf{101} (1994), no.~5, 443--449.

\bibitem{Gathen2003}
Joachim von~zur Gathen and J\"urgen Gerhard, \emph{{Modern Computer Algebra}},
  2nd ed., Cambridge University Press, Cambridge, United Kingdom, 2003.

\bibitem{lcw2008}
Lawrence~C. Washington, \emph{{Elliptic Curves: Number Theory and
  Cryptography}}, 2nd ed., Chapman \& Hall/CRC, 2008.

\bibitem{Williams1998}
Hugh~C. Williams, \emph{{\'{E}douard Lucas and Primality Testing}}, Canadian
  Mathematical Society Series of Monographs and Advanced Texts, vol.~22,
  Wiley-Interscience, 1998.

\end{thebibliography}
%%%%%%%%%%%%%%%%%%%%%%%%%%%%%%%%%%%%%%%%%%%%%%%%%%%%%%%%%%%%%%%%%%%%%%%%%%%%%%
\section*{Appendix: Computing roots of unity (by L. Washington)}
Let $q=2\cdot 3^n+1$ be prime. We show how to construct a $q$th root of
unity mod $p$ (where $p$ is some prime)
in polynomial time in $\log p$ for a fixed $q$.

There are several such primes. The values of $n\le 6000$ are 1, 2, 4, 5,
6, 9, 16, 17, 30, 54, 57, 60, 65, 132, 180, 320, 696, 782, 822, 897,
1252, 1454, 4217, 5480 corresponding to the primes $q=7$, 19, 163,
$\dots$. It is reasonable to conjecture that there are
infinitely many such $q$ (this is similar to the conjecture that there are
infinitely many Mersenne primes).

Let $\zeta_q$ be a primitive $q$th root of unity and let $\rho$ be a
primitive cube root of unity. Let $G$ be the Galois group
of $\Q(\zeta_q,\rho)/\Q(\rho,\sqrt{-q})$. Then
$G$ is cyclic of order $(q-1)/2=3^n$. Let $\sigma$ be a generator and let
\[
\sigma_k=\sigma^{3^{n-k}}.
\]
Then $\sigma_k$ generates a subgroup of $G$ of order $3^k$.
The fixed field $K_k$ of $\sigma_k$ is of degree $3^{n-k}$ over
$\Q(\rho,\sqrt{-q})$.

We want to obtain an expression for a $q$th root of unity
that involves only $\sqrt{-q}$ and taking cube roots.
The basic idea is the following.
Suppose we want to compute $r\in K_m$. Let $r_1=r,\; r_2=\sigma_{m+1}(r),\;
r_3=\sigma_{m+1}^2(r)$ be
the Galois conjugates of $r$ over $K_{m+1}$. Let
\[
f=r_1+r_2+r_3,\quad g=r_1+\rho r_2 + \rho^2 r_3, \quad h= r_1+\rho^2 r_2 +
\rho r_3.
\]
Then $\sigma_{m+1}(g)=\rho^2 g$, so $g^3$ is fixed by $\sigma_{m+1}$ and
therefore lies
in $K_{m+1}$. Similarly, $f^3, h^3\in K_{m+1}$.
If we can determine the values of
$f^3, g^3, h^3$, and if we can compute their cube roots, then we know
$f, g, h$ up to cube roots of unity. So, let's assume that we know
$f, g, h$. Then $r_1=(f+g+h)/3$, $r_2=(f+\rho^2 g+\rho h)/3$, $r_3=(f+\rho
g+\rho^2 h)/3$,
so we recover $r_1, r_2, r_3$.

Start with $r=\zeta_q$.
We will actually use the procedure for $r$ and its Galois conjugates
$\sigma^3(r)$, $\sigma^6(r)$, $\sigma^9(r)$, $\dots$.
The above reduces the computation of $\zeta_q$ and its Galois conjugates
to finding the cube roots of certain
elements of $K_1$. In fact, these elements of $K_1$ are $f^3, g^3, h^3$
and their Galois conjugates over $K_n=\Q(\sqrt{-q}, \rho)$.
We then reduce the computation of these elements to finding the cube roots
of certain
elements of $K_2$ and their conjugates. Continuing in this manner, we
eventually
reduce the problem to computing cube roots of elements of $K_n$. Note that
each time
that we formed a sum $g$, we also formed a sum $h$. These are conjugate
via the automorphism
that sends $\rho$ to $\rho^2$ and fixes $\zeta_q$. Therefore, the elements
of $K_n$ that we obtain
are in pairs $z_1, z_2$ that are conjugate over $\Q(\sqrt{-q})$.
Both $z_1+z_2$ and
$(z_1-z_2)/\sqrt{-3}$ are fixed by Gal($K_n/\Q(\sqrt{-q}))$, so
they lie in
$\Q(\sqrt{-q})$. The real and imaginary parts are rational numbers,
and it is easy to
bound the denominators. Therefore, we can recognize these as rational numbers
by floating point computations. Working back through the preceding and
taking the necessary cube roots,
we obtain an expression for $\zeta_q$.

The expression obtained for $\zeta_q$ can be reduced mod $p$.
There will be some ambiguity caused by the cube roots
being determined only up to powers of $\rho$, so we obtain a finite list of
possibilities of $\zeta_q$. Taking their $q$th powers identifies a
primitive $q$th root of unity.

The above is best understood via an example. Let $q=19$. The Galois group
$G$ is generated by
$\sigma$, which maps $\zeta_{19}$ to $\zeta_{19}^4$. Also,
$\sigma_1=\sigma^3$ maps $\zeta_{19}$
to $\zeta_{19}^7$.
Form
\[
f_0=\zeta_{19}+\zeta_{19}^7+\zeta_{19}^{49}.
\]
The Galois conjugates are $f_0, \sigma(f_0), \sigma^2(f_0)$.

It is classical, and easily verified numerically, that
\[
f_0+\sigma(f_0)+\sigma^2(f_0)=\frac{-1+\sqrt{-19}}2.
\]

Define
\begin{eqnarray*}
x_0&=&(f_0+\rho \sigma(f_0)+\rho^2\sigma^2(f_0))^3\\
x_1&=&(f_0+\rho^2 \sigma(f_0)+\rho\sigma^2(f_0))^3.
\end{eqnarray*}
Then $\sigma$ fixes $x_0$ and $x_1$, so they lie in $\mathbb
Q(\sqrt{-19},\rho)$.
Moreover, the map that switches $\rho$ and $\rho^2$ and fixes $\zeta_{19}$
switches
$x_0$ and $x_1$. Therefore, $x_0+x_1$ and $(x_0-x_1)/\sqrt{-3}$ are in
$\Q(\sqrt{-19})$.
Numerical computation shows that
\begin{eqnarray*}
x_0+x_1&=& \frac 12 (19-17\sqrt{-19})\\
\frac{x_0-x_1}{\sqrt{-3}}&=&\frac 12 (-57-9\sqrt{-19}).
\end{eqnarray*}
(Note that these numbers are algebraic integers, so rounding the results
of a floating point computation yields exact answers.)
Therefore,
\begin{eqnarray*}
x_0&=&\frac14(19-17\sqrt{-19}-57\sqrt{-3}+9\sqrt{57})\\
x_1&=&\frac 14(19-17\sqrt{-19}+57\sqrt{-3}-9\sqrt{57}),
\end{eqnarray*}
where $\sqrt{57}=\sqrt{-3}\sqrt{-19}$.

Therefore, since $1+\rho+\rho^2=0$, we obtain
\[
f_0=\frac 13(\frac{-1+\sqrt{-19}}{2}+x_0^{1/3}+x_1^{1/3}),
\]
with an appropriate choice of cube roots of $x_0$ and $x_1$.

Define
\begin{eqnarray*}
f_1&=&(\zeta_{19}+\rho\sigma_1(\zeta_{19})+\rho^2\sigma_1^2(\zeta_{19}))^3\\
f_2&=&(\zeta_{19}+\rho^2\sigma_1(\zeta_{19})+\rho\sigma_1^2(\zeta_{19}))^3.
\end{eqnarray*}
Then $f_1$ and $f_2$ are fixed by $\sigma_1$, hence lie in $K_1$.
Let
\begin{eqnarray*}
y_1&=&f_1+\sigma(f_1)+\sigma^2(f_1)\\
y_2&=&f_2+\sigma(f_2)+\sigma^2(f_2).
\end{eqnarray*}
Then $y_1$ and $y_2$ lie in $\Q(\sqrt{-19},\rho)$. Numerical
computation
yields
\begin{eqnarray*}
y_1+y_2&=&38-\sqrt{-19}\\
\frac{y_1-y_2}{\sqrt{-3}} &=& 3\sqrt{-19},
\end{eqnarray*}
hence
\begin{eqnarray*}
y_1&=&\frac12 (38-\sqrt{-19}-3\sqrt{57})\\
y_2&=&\frac12 (38-\sqrt{-19}+3\sqrt{57}).
\end{eqnarray*}

Let
\begin{eqnarray*}
x_2&=&(f_1+\rho\sigma(f_1)+\rho^2\sigma^2(f_1))^3\\
x_3&=&(f_1+\rho^2\sigma(f_1)+\rho\sigma^2(f_1))^3\\
x_4&=&(f_2+\rho\sigma(f_2)+\rho^2\sigma^2(f_2))^3\\
x_5&=&(f_2+\rho^2\sigma(f_2)+\rho\sigma^2(f_2))^3.
\end{eqnarray*}
Then
\begin{eqnarray*}
x_2+x_5&=&\frac12(-1007+4373\sqrt{-19})\\
\frac{x_2-x_5}{\sqrt{-3}}&=&\frac12 (-10659-99\sqrt{-19})\\
x_3+x_4&=&1292-1121\sqrt{-19}\\
\frac{x_3-x_4}{\sqrt{-3}}&=&2850+171\sqrt{-19}.
\end{eqnarray*}
Solving yields
\begin{eqnarray*}
x_2&=&\frac14(-1007+4373\sqrt{-19}-10659\sqrt{-3}+99\sqrt{57}) \\
x_3&=&\frac12(1292-1121\sqrt{-19}+2850\sqrt{-3}-171\sqrt{57}) \\
x_4&=&\frac12(1292-1121\sqrt{-19}-2850\sqrt{-3}+171\sqrt{57}) \\
x_5&=&\frac14(-1007+4373\sqrt{-19}+10659\sqrt{-3}-99\sqrt{57}).
\end{eqnarray*}
Again, since $1+\rho+\rho^2=0$, we have
\begin{eqnarray*}
f_1&=&\frac13(x_2^{1/3}+x_3^{1/3}+y_1)\\
f_2&=&\frac13(x_4^{1/3}+x_5^{1/3}+y_2),
\end{eqnarray*}
with an appropriate choice of cube roots.
The search for the appropriate cube roots can
be shortened, for example, by using the fact that
$x_2^{1/3}x_5^{1/3}$ is fixed by $\sigma$ and is unchanged
under the automorphism that maps $\rho$ to $\rho^2$
and which fixes $\zeta_{19}$. It therefore lies
in $\Q(\sqrt{-19})$. Numerical computations show that
$x_2^{1/3}x_5^{1/3}=-114-4\sqrt{-19}$. Therefore, the choice of
cube root for one of $x_2^{1/3}$ and $x_5^{1/3}$ determines the other.

Putting all of the above together, we obtain
\[
\zeta_{19}=\frac13(f_0+f_1^{1/3}+f_2^{1/3})
\]
with an appropriate choice of cube roots.

Schoof's square root algorithm allows us to calculate
$\sqrt{-3}$ and $\sqrt{-19}$
in time polynomial in $\log p$. If taking cube roots mod $p$ is easy
(for example, if $p\equiv 4,\, 7\pmod 9$), then the above quickly calculates
several possibilities for $\zeta_{19}$, corresponding to the choices
of cube roots. Each possibility can be tested to determine whether
or not it is a primitive $19$th root of unity.
This will yield the desired $\zeta_{19}$ in time polynomial in $\log p$.
\end{document}